\documentclass[11pt]{article}
\usepackage{amsfonts}
\usepackage{hyperref}
\usepackage{epsfig}
\usepackage{graphicx}

\newtheorem{theorem}{Theorem}[section]
\newtheorem{definition}[theorem]{Definition}
\newtheorem{proposition}[theorem]{Proposition}
\newtheorem{lemma}[theorem]{Lemma}
\newtheorem{claim}[theorem]{Claim}

\newtheorem{notation}[theorem]{Notation}

\newenvironment{proof}[1][Proof]{ \noindent \textbf{#1: }}{$\Box$
\bigskip}

\begin{document}

\title{Characterization of Co-Blockers for Simple Perfect
Matchings in a Convex Geometric Graph}

\author{
Chaya Keller and Micha A. Perles \\
Einstein Institute of Mathematics, Hebrew University\\
Jerusalem 91904, Israel\\
}

\maketitle

\begin{abstract}
Consider the complete convex geometric graph on $2m$ vertices,
$CGG(2m)$, i.e., the set of all boundary edges and diagonals of a
planar convex $2m$-gon $P$. In~\cite{Maamar-Master}, the smallest
sets of edges that meet all the simple perfect matchings (SPMs) in
$CGG(2m)$ (called ``blockers'') are characterized, and it is shown
that all these sets are caterpillar graphs with a special
structure, and that their total number is $m \cdot 2^{m-1}$. In
this paper we characterize the co-blockers for SPMs in $CGG(2m)$,
that is, the smallest sets of edges that meet all the blockers. We
show that the co-blockers are exactly those perfect matchings $M$
in $CGG(2m)$ where all edges are of odd order, and two edges of
$M$ that emanate from two adjacent vertices of $P$ never cross. In
particular, while the number of SPMs and the number of blockers
grow exponentially with $m$, the number of co-blockers grows
super-exponentially.
\end{abstract}

\section{Introduction}

In this paper we consider convex geometric graphs (i.e., graphs
whose vertices are points in convex position in the plane, and
whose edges are segments connecting pairs of vertices), and in
particular, the complete convex geometric graph on $2m$ vertices,
denoted by $CGG(2m)$.
\begin{definition}
A simple perfect matching (SPM) in $CGG(2m)$ is a set of $m$
pairwise disjoint edges (i.e., edges that do not intersect, not
even in an interior point).
\end{definition}

In~\cite{Maamar-Master}, Keller and Perles give a complete
characterization of the smallest sets of edges in $CGG(2m)$ that
meet all the SPMs, called {\it blockers}. It turns out that all
the blockers are simple trees of size $m$ admitting a special
structure called {\it caterpillar graphs}~\cite{Caterpillar0,Caterpillar1},
and that their number is $m \cdot 2^{m-1}$.

Following the result of~\cite{Maamar-Master}, one may consider a
sequence $\{A_n\}_{n=0}^{\infty}$, defined inductively as follows.
$A_0$ is the family of all SPMs in $CGG(2m)$. Given $A_k$, define
$A_{k+1}$ as the family of all smallest sets of edges in $CGG(2m)$
that meet all of the elements of $A_k$. In particular, $A_1$ is
the family of all blockers, characterized in~\cite{Maamar-Master}.

A standard argument shows that $A_3=A_1$, and thus $A_k=A_{k-2}$
for all $k \geq 3$. Thus, the only unknown element of the sequence
is $A_2$, i.e., the family of all smallest sets of edges of
$CGG(2m)$ that meet all blockers, called in the sequel {\em
co-blockers}. It is easy to show (see Section~\ref{sec:main}) that
the size of any co-blocker is at least $m$, and on the other hand,
any SPM meets every blocker by the definition of a blocker, and
thus is a co-blocker. Therefore, the size of the co-blockers is
$m$ (like the size of the SPMs and of the blockers).

\begin{figure}[tb]
\begin{center}
\scalebox{0.8}{
\includegraphics{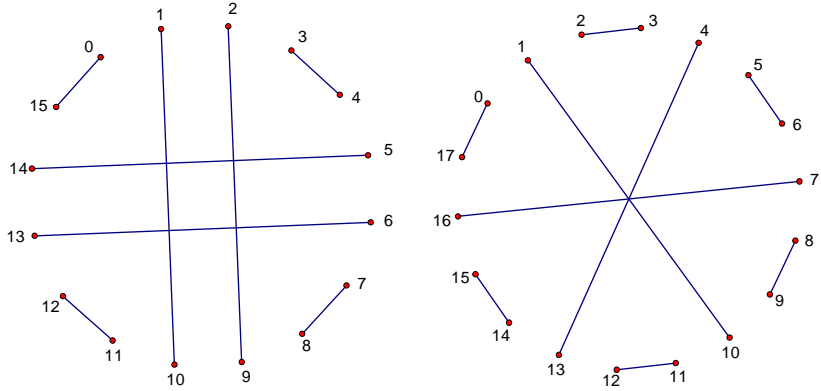}
} \caption{Two small co-blockers that are not SPMs.}
\label{Fig1}
\end{center}
\end{figure}

In this paper we give a complete characterization of the family of
co-blockers:
\begin{theorem}\label{Thm:Main}
For any $m \in \mathbb{N}$, the set of co-blockers in $CGG(2m)$ is
the set of all perfect matchings in $CGG(2m)$, such that:
\begin{itemize}
\item All the edges of the matching have odd order (see
Section~\ref{sec:preliminaries} for a formal definition of the
order of an edge in $CGG(2m)$).

\item Two edges $[a,b]$ and $[a',c]$ of the matching whose
end-points $a,a'$ form a boundary edge of $CGG(2m)$ never cross.
\end{itemize}
\end{theorem}
Two examples of small co-blockers that are not SPMs are given in
Figure~\ref{Fig1}.

\medskip

\noindent The theorem allows us to give lower and upper bounds on the number
of co-blockers:
\begin{proposition}\label{Prop:Number0}
Denote the set of co-blockers in $CGG(2m)$ by $A_2(m)$. Then for
all $m \in \mathbb{N}$,
\[
\lfloor m/3 \rfloor ! \leq |A_2(m)| \leq m!.
\]
\end{proposition}

It is known that both the number of SPMs and the number of
blockers grow only exponentially with $m$: it is easy to show that the
number of SPMs is the Catalan number
$C_m=\frac{1}{m+1}{{2m}\choose{m}}$, and it is shown
in~\cite{Maamar-Master} that the number of blockers is $m \cdot
2^{m-1}$. Thus, Proposition~\ref{Prop:Number0} shows that the
number of co-blockers is significantly larger than the numbers of
SPMs and blockers.

The rest of this paper is organized as follows: In
Section~\ref{sec:preliminaries} we introduce some basic
definitions and recall the properties of SPMs and blockers that
are used in our proof. In Section~\ref{sec:main} we present the
proof of Theorem~\ref{Thm:Main}. Finally, in
Section~\ref{sec:number} we prove Proposition~\ref{Prop:Number0}.

\section{Preliminaries}
\label{sec:preliminaries}

In this section we introduce several basic definitions, and recall
some properties of SPMs and of blockers presented
in~\cite{Maamar-Master}, which are used in the proof of
Theorem~\ref{Thm:Main}.

\subsection{Definitions and Notations}
\label{sec:sub:definitions}

Throughout this paper, we use the following definitions and
notations.

\begin{notation}
The set of vertices of $CGG(2m)$ is denoted by $V$, and is
realized in the plane as the set of vertices of a convex $2m$-gon
$P$. The vertices are labelled cyclically from $0$ to $2m-1$.
\end{notation}

\begin{definition}
The \textbf{order} of an edge $[i,j]$ is $\min(|j-i|,2m-|j-i|)$.
The boundary edges of $P$ are, of course, of order $1$. We call
the non-boundary edges, i.e., the edges that are diagonals of $P$,
\textbf{interior edges}.
\end{definition}

\begin{definition}
The \textbf{direction} of an edge in $CGG(2m)$ is the sum (modulo
$2m$) of the labels of its endpoints. That is, if $e=[i,j]$, then
its direction is:
\[
Dir(e)= i+j (\mbox{ mod } 2m) = \left\lbrace
  \begin{array}{c l}
    i+j, & i+j<2m\\
    i+j-2m, & i+j \geq 2m.
  \end{array}
\right.
\]
Two edges $e,e'$ are \textbf{parallel} if
$Dir(e)=Dir(e')$.\footnote{Note that if $P$ is regular, an
equivalent definition is that $e,e'$ are parallel as straight line
segments in the plane.}
\end{definition}

\begin{definition}
Two edges $e,e'$ of $CGG(2m)$ are called \textbf{neighbors} if (at
least) one endpoint of $e$ is adjacent to (at least) one
endpoint of $e'$ on the boundary of $P$.
\end{definition}

\begin{definition}
A perfect matching $M$ of $CGG(2m)$ is called \textbf{semi-simple}
if:
\begin{itemize}
\item All the edges of $M$ are of odd order, and

\item $M$ does not contain a pair of crossing neighbors.
\end{itemize}
\end{definition}

\subsection{Caterpillar Trees and the Structure of Blockers}
\label{sec:sub:blockers}

\begin{definition}
A tree $T$ is a \textbf{caterpillar} (or a fishbone) if the
derived graph $T'$ (i.e., the graph obtained from $T$ by removing
all leaves and their incident edges) is a path (or is empty). A
geometric caterpillar is \textbf{simple} if it does not contain a
pair of crossing edges. A longest path in a caterpillar $T$ is
called a \textbf{spine} of $T$. Given a spine of $T$, the edges of
$T$ that have one endpoint interior to the spine and the other
endpoint exterior to the spine are called \textbf{legs} of $T$.
\end{definition}

In~\cite{Maamar-Master}, the blockers in $CGG(2m)$ are fully
characterized by the following theorem:
\begin{theorem}\label{Thm:Blockers}
Any blocker of $CGG(2m)$ is a simple caterpillar graph whose spine
lies on the boundary of $P$ and is of length $t \geq 2$. If the
spine ``starts'' with the vertex $0$ and the edge $[0,1]$, then
the edges of the blocker are:
\begin{equation}\label{Eq2.1}
\{[i-1,i]:1 \leq i \leq t\} \cup
\{[t+j-1-\epsilon_{t+j},t+j+\epsilon_{t+j}]:1 \leq j \leq m-t\},
\end{equation}
where the $\epsilon_i$ are natural numbers satisfying $1 \leq
\epsilon_{t+1} < \epsilon_{t+2} < \ldots < \epsilon_m \leq m-2$.

\medskip

\noindent Conversely, any set of $m$ edges of the described form is a
blocker in $CGG(2m)$.
\end{theorem}

\begin{figure}[tb]
\begin{center}
\scalebox{0.6}{
\includegraphics{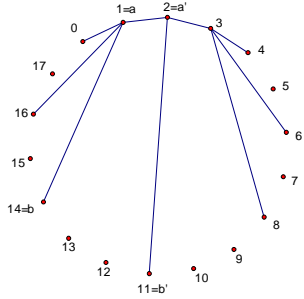}
} \caption{A blocker in $CGG(18)$.}
\label{Figure2}
\end{center}
\end{figure}

The contents of Formula~(\ref{Eq2.1}) can be described as follows:
\begin{enumerate}
\item For each pair $e,e'$ of opposite boundary edges of $P$, the
blocker contains exactly one edge parallel to $e$ and $e'$.

\item Each leg of the caterpillar connects a vertex $a$ interior
to the spine to a vertex $b$ exterior to the spine.

\item If $[a,b]$ and $[a',b']$ are two distinct legs of the
caterpillar (where $a,a'$ are on the spine and $b,b'$ are not),
then the distance between $b$ and $b'$ along the complement of the
spine is larger than the distance between $a$ and $a'$ (within the
spine). In other words, if the spine starts with the vertex $0$
and $a<a'$, then
\begin{equation}\label{NewFormula1}
b-b'>a'-a
\end{equation}
(see Figure~\ref{Figure2}).
\end{enumerate}

An example of a blocker in $CGG(18)$ is depicted in Figure~\ref{Figure2}.

\medskip

\noindent In our proof we also use the following simple claim:
\begin{claim}\label{Cor:Blockers}
In $CGG(2m)$, the set of all edges of odd order emanating from a single
vertex is a blocker. This blocker is called ``a star blocker''.
\end{claim}

The star blockers correspond to the case $t=2$ in
Theorem~\ref{Thm:Blockers}. The other extreme value $t=m$ yields
blockers that are just one half of the boundary circuit of $P$.


\section{Proof of Theorem~\ref{Thm:Main}}
\label{sec:main}

In this section we present the proof of our main theorem. We start
by observing a simple necessary condition for co-blockers.
As we shall see later, this condition is not so far from being
sufficient.
\begin{lemma}\label{Lemma:Basic}
Let $C$ be a co-blocker in $CGG(2m)$. Then $C$ is a perfect
matching, and all the edges of $C$ are of odd order.
\end{lemma}

\begin{proof}
First, note that if there exists a vertex $x$ that is not
contained in any edge of $C$, then $C$ does not meet the star
blocker emanating from $x$, contradicting the assumption that $C$
is a co-blocker. Thus, any vertex $x \in V$ is contained in an
edge of $C$, and since $C$ has only $m$ edges (as noted in
the introduction), this implies that $C$ is a perfect matching.

Second, suppose on the contrary that $C$ contains an edge
$e=[x,y]$ of even order. Since $e$ is the only edge of $C$ that
emanates from $x$, we again find that $C$ does not meet the star
blocker emanating from $x$, contradicting the assumption. Thus,
all edges of $C$ are of odd order.
\end{proof}

We proceed by observing a property of semi-simple perfect
matchings which will be crucial in our analysis:
\begin{lemma}\label{Obs:Key}
Let $M$ be a semi-simple perfect matching in $CGG(2m)$. Then the
following holds:
\begin{itemize}
\item If $e$ is an interior edge in $M$, then $M$ contains a
boundary edge in each of the two open half-planes determined by
the straight line \emph{\textrm{aff}}$(e)$.

\item If $e_1,e_2$ are two crossing edges of $M$ (i.e., edges
which intersect in an interior point), then $M$ contains a
boundary edge in each of the four open quadrants determined by
\emph{\textrm{aff}}$(e_1)$ and \emph{\textrm{aff}}$(e_2)$.
\end{itemize}
\end{lemma}

\begin{figure}[tb]
\begin{center}
\scalebox{0.8}{
\includegraphics{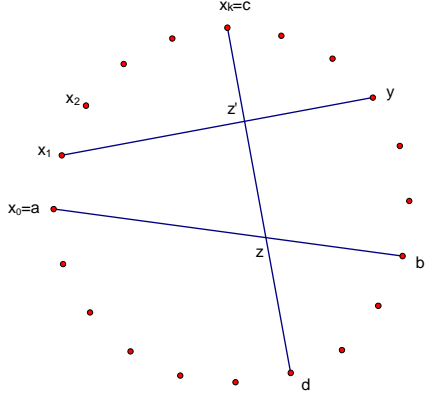}
} \caption{Illustration of the proof of Lemma~\ref{Obs:Key}.}
\label{Figure3}
\end{center}
\end{figure}

\begin{proof}
We begin with the first claim. Let $e=[a,b]$ and let $H$ be one of
the half-planes determined by \textrm{aff}$(e)$. $H$ meets the
boundary of $P$ in a polygonal arc $\langle x_0,x_1,\ldots,x_k
\rangle$, where $x_0=a$ and $x_k=b$. Consider the set of all edges
of $M$ both of whose endpoints are in $\{x_0,x_1,\ldots,x_k\}$,
like $e$. Among those edges, choose an edge $e'=[x_i,x_j]$ ($i<j$)
that minimizes the difference $j-i$. We claim that $e'$ is a
boundary edge.

Indeed, if $e'$ is not a boundary edge, then $x_{i+1}$ is an
internal vertex of the polygonal arc $\langle
x_i,x_{i+1},\ldots,x_j \rangle$. Let $e''$ be the edge of $M$ that
contains $x_{i+1}$. By the minimality of $e'$, the other endpoint
of $e''$ cannot be in $\{x_i,x_{i+1},\ldots,x_j\}$, and thus, $e'$
and $e''$ are crossing neighbors in $M$, contradicting the
assumption that $M$ is semi-simple. Hence, $e'$ is indeed a
boundary edge, as asserted.

\medskip

\noindent Now we proceed to the second claim. Let $e_1=[a,b],
e_2=[c,d]$, and $e_1 \cap e_2 = \{z\}$. Let $Q$ be the quadrant
determined by the rays $\overrightarrow{za}$ and
$\overrightarrow{zc}$. $Q$ meets the boundary of $P$ in a
polygonal arc $\langle x_0,x_1,\ldots,x_k \rangle$, where $x_0=a$
and $x_k=c$. We proceed by induction on $k$. The case $k=1$ is
impossible, since otherwise $e_1$ and $e_2$ are crossing
neighbors, which contradicts the assumption that $M$ is
semi-simple.

Thus, we may assume that $k \geq 2$, and, in particular, that
$x_1$ is an internal vertex of the polygonal arc $\langle
x_0,x_{1},\ldots,x_k \rangle$. Let $e'=[x_1,y]$ be the edge of $M$
that contains $x_{1}$. We consider four cases:
\begin{itemize}
\item If $y$ is in $\{x_0,x_{1},\ldots,x_k\}$, then by the first
claim, $M$ contains a boundary edge $e''$ both of whose endpoints
are in $\{x_{1},x_2,\ldots,x_{k-1}\}$, hence $e'' \subset$
{\textrm{int}}$(Q)$.

\item If $y$ is on the boundary of $P$ strictly between $x_k$
    and $b$, then the edge $[x_1,y]$ crosses the edge $[c,d]$ at
    some point $z' \in$ \textrm{int}$(P)$ (see
    Figure~\ref{Figure3}). The quadrant $Q'$ determined by the
    rays $\overrightarrow{z'x}$ and $\overrightarrow{z'c}$ meets
    the boundary of $P$ in a shorter polygonal arc $\langle
    x_1,x_{2},\ldots,x_k \rangle$. Thus, by the induction
    hypothesis, $M$ contains a boundary edge in
    \textrm{int}$(Q')$, and that edge is (of course) contained
    in \textrm{int}$(Q)$.

\item If $y=b$, then $M$ contains two edges emanating from the
same vertex, contradicting the assumption that $M$ is a perfect
matching.

\item If $y$ is not one of the above, then the edges $[a,b]$ and
$[x_1,y]$ are crossing neighbors, contradicting the assumption
that $M$ is semi-simple.
\end{itemize}

This completes the proof.
\end{proof}

Now we are ready to state our main theorem.
\begin{theorem}\label{Thm:Main1}
Let $M$ be a set of $m$ edges in $CGG(2m)$. Then $M$ is a
co-blocker if and only if $M$ is a semi-simple perfect matching.
\end{theorem}

\begin{figure}[tb]
\begin{center}
\scalebox{0.8}{
\includegraphics{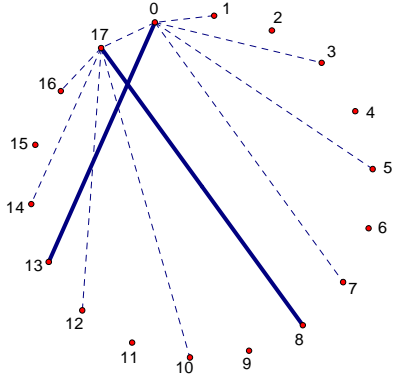}
} \caption{A blocker that misses a perfect matching with crossing
neighbors.} \label{Figure4}
\end{center}
\end{figure}

\begin{proof}
\textbf{Necessity:} Assume that $M$ is a co-blocker. By
Lemma~\ref{Lemma:Basic}, $M$ is a perfect matching and all its
edges are of odd order. Suppose on the contrary that $M$ is not
semi-simple, and thus w.l.o.g. contains the crossing neighbors
$[0,2l-1]$ and $[2k,2m-1]$, where $0<2k<2l-1<2m-1$. Then the
blocker $B$ whose spine is $\langle 2m-2,2m-1,0,1 \rangle$ and
whose legs are $[0,2j-1]$ for all $2 \leq j < l$ and $[2j,2m-1]$
for all $l \leq j < m-1$ does not meet $M$, contradicting the
assumption that $M$ is a co-blocker. The blocker $B$ is depicted
in Figure~\ref{Figure4}.

\medskip

\textbf{Sufficiency:} Assume $M$ is a semi-simple perfect
matching, and suppose on the contrary that $M$ misses some blocker
$B$. Assume, without loss of generality, that the spine of $B$ is
$\langle 0,1,2,\ldots,t \rangle$, where $2 \leq t \leq m$. (Note
that by Theorem~\ref{Thm:Blockers}, the blocker is a caterpillar
whose spine lies on the boundary of $P$.) For $i=1,2,\ldots,t-1$,
denote by $e_i$ the (unique) edge of $M$ that emanates from $i$,
and denote its other endpoint by $y_i$. We claim that the edges
$e_1,e_2,\ldots,e_{t-1}$ satisfy the following:
\begin{itemize}
\item For any $1 \leq i \leq t-1$, we have $y_i \in \{t+1,t+2,
\ldots,2m-1\}$. In particular, the $t-1$ edges
$e_1,e_2,\ldots,e_{t-1}$ are distinct.

\item For any pair $i,j \in \{1,2,\ldots,t-1\}$, the edges $e_i$
and $e_j$ do not cross.
\end{itemize}

The first claim follows from the first claim of
Lemma~\ref{Obs:Key}. Indeed, if $y_i \in \{0,1,\ldots,t\}$, then
by the lemma, $M$ contains a boundary edge in the polygonal path
between $i$ and $y_i$, and by assumption, this edge is in the
spine of $B$, which contradicts the assumption that $M$ misses
$B$. The second claim follows similarly from the second claim of
Lemma~\ref{Obs:Key}.

The second claim implies that if $1 \leq i<i+1 \leq t-1$, then
$y_{i+1}<y_i$, and therefore, $i+y_i \geq (i+1)+y_{i+1}$. Thus,
the function $g:i \mapsto i+y_i$ is monotone non-increasing in $i$
for $1 \leq i \leq t-1$. We can also bound the range of this function,
namely,
\[
2m-1=1+(2m-2) \geq g(1) \geq g(i) \geq g(t-1) \geq (t-1)+(t+2) =
2t+1,
\]
which implies that for all $1 \leq i \leq t-1$,
\[
Dir(e_i) = i+y_i (\mbox{ mod } 2m) = i+y_i =g(i).
\]
Since, by Theorem~\ref{Thm:Blockers}, the blocker $B$ contains a
unique edge parallel to every edge of odd order in $CGG(2m)$,
there is a unique edge $f_i$ of $B$ parallel to $e_i$. This edge
cannot lie on the spine of $B$, since for the edges on the spine
of $B$, the direction takes the values $1,3,\ldots,2t-1$, and thus
they are not parallel to the edges $e_i$. Hence, $f_i$ is a leg of
$B$, which can be represented as $f_i=[r_i,q_i]$, with $1 \leq r_i
\leq t-1$, and $t+1 \leq q_i \leq 2m-1$.

Note that we have $Dir(f_i)=r_i+q_i$. (The other option,
$Dir(f_i)=r_i+q_i-2m$, would yield $Dir(f_i) \leq t-2$, whereas
$Dir(f_i)=Dir(e_i) \geq 2t+1$.) Thus, by the monotonicity of
$Dir(e_i)$, we have:
\begin{equation}\label{NewFormula2}
r_{i+1}+q_{i+1} = Dir(f_{i+1}) = Dir(e_{i+1}) \leq Dir(e_i) =
Dir(f_i) = r_i+q_i,
\end{equation}
for $i=1,2,\ldots,t-2$.

\medskip

Now we are ready to reach the contradiction. If $r_i=i$ for some
$i$, then $f_i=e_i$, contrary to the assumption that $M \cap B =
\emptyset$. Thus, $r_1>1$ and $r_{t-1}<t-1$. Hence, there are two
consecutive indices $1 \leq i<i+1 \leq t-1$ with $r_i>i$ and
$r_{i+1}<i+1$. It follows that $r_{i+1} < r_i$, and therefore, by
Theorem~\ref{Thm:Blockers}, we must have
$q_{i+1}-q_{i}>r_{i}-r_{i+1}$ (see Equation~(\ref{NewFormula1})),
which contradicts Equation~(\ref{NewFormula2}). This completes the
proof.
\end{proof}

\section{The Number of Co-Blockers}
\label{sec:number}

The characterization of the co-blockers given in
Theorem~\ref{Thm:Main} allows us to find upper and lower bounds on
the number of co-blockers in $CGG(2m)$, as a function of $m$.

\begin{proposition}
The number of co-blockers satisfies
\begin{equation}\label{Eq:Number1}
\lfloor m/3 \rfloor ! \leq |C(m)| \leq m!.
\end{equation}
\end{proposition}
\begin{proof}
The right inequality in~(\ref{Eq:Number1}) follows immediately
from Lemma~\ref{Lemma:Basic}, since by the lemma, all the
co-blockers are perfect matchings whose edges are of odd order.
These matchings can be viewed as bijections from the set of
vertices of odd index to the set of vertices of even index, and
their number is $m!$.

\medskip

In order to prove the left inequality in~(\ref{Eq:Number1}), we
consider perfect matchings of a special type. For the sake of
simplicity we assume first that $m=3k$, and denote the vertices of
$CGG(2m)$ by $0,1,2,\ldots,2m-1$. We consider only perfect
matchings that contain all the boundary edges
\[
[1,2],[4,5],\ldots,[6k-5,6k-4],[6k-2,6k-1],
\]
i.e., all the boundary edges whose vertices are congruent to $1$
and $2$ modulo 3. We claim that any perfect matching of this class
whose edges are all of odd order is a co-blocker. Indeed, by
Theorem~\ref{Thm:Main}, such a perfect matching is not a
co-blocker only if it contains two edges whose endpoints are
consecutive vertices on the boundary which intersect in an
interior point. However, amongst any two consecutive vertices on
the boundary there is a vertex whose index modulo 3 equals to $1$
or $2$, and thus the edge of the matching containing that vertex
is a boundary edge and cannot cross any other edge.

The number of perfect matchings of this class is $k!$, since any
vertex whose index equals $0$ modulo $6$ can be connected by an
edge to any vertex whose index equals $3$ modulo $6$, and all the
other vertices are already contained in boundary edges.

Thus, $|C(m)| \geq k! = (m/3)!$. Finally, if $m=3k+1$ or $m=3k+2$,
then one may consider perfect matchings of the class described
above, but containing also the boundary edge $[6k,6k+1]$ (for both
$m=3k+1,m=3k+2$), and in addition $[6k+2,6k+3]$ (for $m=3k+2$).
The argument given above in the case $m=3k$ applies also here, and
the number of such perfect matchings is $k!$. Thus, $|C(m)| \geq
\lfloor m/3 \rfloor !$, as asserted.
\end{proof}

\end{document}